\documentclass{svproc}
\usepackage{amsmath,amssymb,latexsym} 
\usepackage{url}

\newtheorem{algorithm}{Algorithm}
\begin{document}
\mainmatter              
\title{Complete solving the quadratic equation mod $2^n$}
\titlerunning{Complete solving the quadratic equation mod $2^n$}  
%
\author{S. M. Dehnavi\inst{1} \and M. R. Mirzaee Shamsabad\inst{2} \and A. Mahmoodi Rishakani\inst{3}}
\authorrunning{S. M. Dehnavi et al.} 
%
\tocauthor{S. M. Dehnavi, M. R. Mirzaee Shamsabad, and A. Mahmoodi Rishakani}
\institute{Department of Mathematical and Computer Sciences, University of Kharazmi, Tehran, Iran,\\
\email{dehnavism@ipm.ir}
\and
Department of Mathematics, Shahid Beheshti University, Tehran, Iran,\\
\email{m\_mirzaee@sbu.ac.ir}
\and
Department of Sciences, Shahid Rajaee Teacher Training University, Tehran, Iran,\\
\email{am.rishakani@srttu.edu}}

\maketitle              

\begin{abstract} 
Quadratic functions have applications in cryptography. In this paper, we investigate the modular quadratic equation 
$$ ax^2+bx+c=0 \quad (mod \,\, 2^n), $$
and provide a complete analysis of it.  More precisely, we determine when this equation has a solution and in the case that it has a solution, we not only determine the number of solutions, but also give the set of solutions in $O(n)$ time. One of the interesting results of our research is that, when this equation has a solution, then the number of solutions is a power of two.


\end{abstract}
\section{Introduction and Main Results}\label{SecInt}
The square mapping is one of the tools which is used in cryptography. For instance, the Rabin cryptosystem~\cite{ST} employs modular quadratic mapping. As another example, in the design of the stream cipher Rabbit~\cite{BOVE}, the square map is used.  A quadratic polynomial modulo $2^{32}$ is used in the AES finalist block cipher RC6~\cite{RIRO}.\newline
Quadratic equation has been solved over various algebraic structures. For example, the quadratic equation over $\mathbb{F}_{2^n}$ is solved in Theorem 3.2.15 of~\cite{GOKL}. Note that an algorithm for finding the solutions of quadratic equation over $\mathbb{F}_{2^n}$ is also given in \cite{PO}. This research is not the first one concerning the quadratic equation mod $2^n$. For instance, \cite{Q} gives the number of solutions of equation~(\ref{QuadEqu}) in spacial cases,  but the authors do not give the set of solutions.
\newline
In this paper, we examine the quadratic equation mod $2^n$. We verify when this equation has a solution and in the case that it has a solution, we give the number of solutions as well as the set of its solutions in $O(n)$ time.
\newline
In section 2, we give the preliminary notations and definitions. Section 3 is devoted to the main theorems of the paper which solve the modular quadratic equation mod $2^n$ completely and presents its number of solutions along with its set of solutions. In section 4, we conclude the paper. 

\section{Notations and Definitions} \label{SecNot}
Throughout this paper $n$, $a$, $b$, and $c$ are natural numbers. We denote the well-known ring of integers mod $2^n $ by $\mathbb{Z}_{2^n}$. For every nonzero element $a \in \mathbb{Z}_{2^n}$, we define $p_2(a)$ as the greatest power of 2 that divides $a$ and the odd part of $a$ or $\dfrac{a}{2^{p_2(a)}}$ is denoted by $o_2(a)$. Here, we define $p_2(0)=n$. 
\newline
The number of elements (cardinal) of a finite set $ A $ is denoted by $ \vert A \vert $. For a function $ f :A \to B $, the preimage of an element $ b \in B $ is defined as $ \{ a \in A \vert b=f(a) \} $ and denoted by $ f^{-1}(b) $. The $i$-th bit of a natural number $x$ is denoted by $[x]_i$.  Also, for an integer $j$, we define $e_j$ as follows
$$ 
e_j= 
\begin{cases} 
0 & \text{j odd}\\ 
1 & \text{j even}\\ 
\end{cases} 
$$ 
Let $(G,*)$ be a group and $\varphi:G \to G$ be a group endomorphism. We denote the kernel of $\varphi$ by $ker(\varphi)$ and the image of $\varphi$ by $Im(\varphi)$. 

\section{Solving the quadratic equation mod $2^n$}\label{SecSol}
	
In this section, we study the modular quadratic equation 
\begin{equation}\label{QuadEqu}
ax^2+bx+c=0 \quad (mod \,\, 2^n),
\end{equation}
and wish to solve it. More precisely, we want to determine:
\newline
a) whether (\ref{QuadEqu}) has a solution.
\newline
b) if it has a solution, then what is the number of its solutions.
\newline
c) the set of solutions.
\newline
In the sequel, we note that $x=0$ is equivalent to $p_2 (x)=n $.
\begin{lemma}\label{abceven}
	Let $a$, $b$, and $c$ be even and $t=min \{ p_2(a), \, p_2(b), \, p_2(c) \}$. Set $\acute{a}={a \over 2^t}$, $\acute{b}={b \over 2^t}$, and $\acute{c}={c \over 2^t}$. Consider the equations~(\ref{QuadEqu}) and 
	\begin{equation}\label{primeQuadEqu}
		\acute{a}x^2+\acute{b}x+\acute{c}=0 \quad (mod \,\, 2^{n-t}).
	\end{equation}
	Len $N_1$ and $N_2$ be the number of solutions of (\ref{QuadEqu}) and (\ref{primeQuadEqu}) respectively. Also, let $\{ x_1,\dots,x_{N_2} \}$ be the set of solutions of (\ref{primeQuadEqu}). Then the set of solutions of (\ref{QuadEqu}) is as follows
	$$\{x_i+r2^{n-t}:\,\, 0\leq r < 2^t, \,\, 1 \leq i \leq N_2 \}.$$
	Further $N_1=2^tN_2$.
	\begin{proof}
		Firstly, fix $1 \leq i \leq N_2$ and $0\leq r < 2^t$. We show that $x_i+r2^{n-t}$ is a solution of (\ref{QuadEqu}):
		$$
		\begin{array}{llll}
			a ( x_i+r2^{n-t} )^2+b ( x_i+r2^{n-t} ) + c\\
			=2^t \acute{a} ( x_i^2+r^22^{2n-2t}+x_ir2^{n-t+1}) + 2^t \acute{b} ( x_i+r2^{n-t}) + 2^t \acute{c}\\
			= 2^t \acute{a} x_i^2 + \acute{a} r^2 2^{2n-t} + \acute{a} x_i r 2^{n+1} + \acute{b} x_i 2^t + \acute{b} r 2^n + 2^t \acute{c}\\
			= 2^t ( \acute{a} x_i^2 + \acute{b} x_i + \acute{c}) = 0 \quad (mod \,\, 2^n).
		\end{array}
		$$
		Conversely, let $x \in \mathbb{Z}_{2^n}$ be a solution of (\ref{QuadEqu}). Then
		$$2^t ( \acute{a} x^2 + \acute{b} x + \acute{c}) = 0 \,\, (mod \,\, 2^n),$$
		so
		$$\acute{a}x^2+\acute{b}x+\acute{c}=0 \,\, (mod \,\, 2^{n-t}),$$
		and one can check that $\chi =x \,\, (mod\,\,2^{n-t})$ is a solution of (\ref{primeQuadEqu}). Thus all of solutions $y$ of (\ref{QuadEqu}) are such that $y = x_i + r2^{n-t}$ for some $1 \leq i \leq N_2$ and $0\leq r < 2^t$.
	\end{proof}
\end{lemma}
\begin{example}
	Consider the equation
	\begin{equation}\label{Exa1}
		4x^2+4x+24=0 \quad (mod \,\, 2^5),
	\end{equation}
	and
	\begin{equation}\label{Exa2}
		x^2+x+6=0 \quad (mod \,\, 2^3).
	\end{equation}
	The set of solutions of (\ref{Exa1}) and (\ref{Exa2}) are $A=\{1,6,9,14,17,22,25,30\}$ and $B=\{1,6\}$, respectively. It is not hard to see that Lemma~\ref{abceven} holds for this example and $|A|=2^2|B|$.
\end{example}
The proof of following lemma is straightforward.
\begin{lemma}\label{NoSol}
	The equation~(\ref{QuadEqu}) has no solutions when $p_2(a)=p_2(b)=p_2(c)=0$ or when $p_2(a)>0$, $p_2(b)>0$, and $p_2(c)=0$.
\end{lemma}	
\begin{lemma}\label{aevenbodd}
	If $p_2(a)>0$ and $p_2(b)=0$, then the equation~(\ref{QuadEqu}) has a unique solution.
	\begin{proof}
		Consider two following cases:\\
		Case I) Let $p_2(a)>0$, $p_2(b)=p_2(c)=0$, $a=2 \acute{a}$, $b=2 \acute{b} +1$, and $c=2 \acute{c}+1$. In this case, any solution $x$ of (\ref{QuadEqu}) is odd; so, we have $x=2 \acute{x} +1$. Thus
		$$ 2 \acute{a}(2 \acute{x}+1)^2 + (2 \acute{b}+1)(2 \acute{x}+1) +2\acute{c}+1=0 \quad (mod \,\, 2^n).$$
		Then
		$$4\acute{a}\acute{x}^2+(4\acute{a}+2\acute{b}+1)\acute{x}+\acute{a}+\acute{b}+\acute{c}+1=0 \quad (mod \,\, 2^{n-1}).$$
		So, if we set $\alpha=4\acute{a}$, $\beta=4\acute{a}+2\acute{b}+1$, and $\gamma=\acute{a}+\acute{b}+\acute{c}+1$, then $[x]_0=1$ and we must solve the equation 
		$$\alpha\acute{x}^2+\beta\acute{x}+\gamma=0 \quad (mod \,\, 2^{n-1}),$$
		such that $p_2(\alpha)>0$ and $p_2(\beta)=0$. Now, we have either $p_2(\gamma)=0$ which is this same case or $p_2(\gamma)>0$ which is Case II bellow.\\
		Case II) Let $p_2(a)>0$, $p_2(c)>0$, $p_2(b)=0$, $a=2\acute{a}$, $b=2\acute{b}+1$, and $c=2\acute{c}$. In this case, $x=2\acute{x}$ and thus
		$$2\acute{a}(2\acute{x})^2+(2\acute{b}+1)(2\acute{x})+2\acute{c}=0 \quad (mod \,\, 2^n).$$
		Then
		$$4\acute{a}\acute{x}^2+(2\acute{b}+1)\acute{x}+\acute{c}=0 \quad (mod \,\, 2^{n-1}).$$
		Put $\alpha=4\acute{a}$, $\beta=2\acute{b}+1$, and $\gamma=\acute{c}$. Then $[x]_0=0$ and we must solve the equation 
		$$\alpha\acute{x}^2+\beta\acute{x}+\gamma=0 \quad (mod \,\, 2^{n-1}),$$
		with $p_2(\alpha)>0$ and $p_2(\beta)=0$. Now, if $p_2(\gamma)=0$, then we transit to Case I and if $p_2(\gamma)>0$, then we transit to this same case.
		Therefore, (\ref{QuadEqu}) has a unique solution.
	\end{proof}
\end{lemma}
The trend of the proof of Lemma~\ref{aevenbodd} justifies the correctness of Algorithm~\ref{Alg}, which computes the solution of (\ref{QuadEqu}) with conditions of Lemma~\ref{aevenbodd} in $O(n)$ time.
\begin{algorithm}\label{Alg}
	\textbf{Solve$(a,b,c,n)$}\newline
	Input: $a,\,b,\,c \in \mathbb{Z}_{2^n}$ with $p_2(a)>0$ and $p_2(b)=0$.\newline
	Output: The solution of (\ref{QuadEqu}) in binary form.\newline
	for $i=0$ to $n-1$ do \newline
	begin\newline
	if $p_2(c)>0$ then
	
	$[x]_i=0$ 
	
	Solve$(2a,b,{c \over 2},n-1)$ \newline
	else
	
	$[x]_i=1$ 
	
	Solve$(2a,2a+b,{a \over 2}+\lfloor{b \over 2}\rfloor + \lfloor {c \over 2}\rfloor +1,n-1)$.
\end{algorithm}
\begin{lemma}\label{aboddceven}
	In the case that $p_2(a)=p_2(b)=0$ and $p_2(c)>0$, the equation (\ref{QuadEqu}) has two solutions.
	\begin{proof}
		Consider the equation $2ay^2+by+{c \over 2}=0 \quad (mod \,\, 2^{n-1}).$ Lemma~\ref{aevenbodd} shows that this equation has a unique solution $\delta\in \mathbb{Z}_{2^{n-1}}$. One can check that $\pi =2 \delta$ is a solution of (\ref{QuadEqu}) in this case. On the other hand, Lemma~\ref{aevenbodd} shows that the equation 
		$$2az^2+(2a+b)z+{a+b+c \over 2}=0 \quad (mod \,\, 2^{n-1}),$$
		has a unique solution $\rho \in  \mathbb{Z}_{2^{n-1}}$. It is straightforward to see that $\varepsilon=2\rho+1$ is a solution of (\ref{QuadEqu}) in $\mathbb{Z}_{2^n}$. Now, we show that (\ref{QuadEqu}) has no other solutions. Suppose that $x$ is a solution of (\ref{QuadEqu}). We have two following cases:\\
		Case I) $x=2\acute{x}$, then
		$$4a\acute{x}^2+2b\acute{x}+c=0 \quad (mod \,\, 2^n),$$
		so
		$$2a\acute{x}^2+b\acute{x}+{c \over 2} =0 \quad (mod \,\, 2^{n-1}),$$
		which is not a new solution.\\
		Case II) $x=2\acute{x}+1$, then
		$$4a\acute{x}^2+(4a+2b)\acute{x}+a+b+c=0 \quad (mod \,\, 2^n),$$
		so
		$$2a\acute{x}^2+(2a+b)\acute{x}+{a+b+c \over 2} =0 \quad (mod \,\, 2^{n-1}),$$
		which is not a new solution.
	\end{proof}
\end{lemma}
The proof of next lemma is straightforward.
\begin{lemma}\label{aoddmod8}
	If $a$ is an odd element in $\mathbb{Z}_{2^n}$. Then $a^2=1 \,\, (mod \,\, 8)$. 
\end{lemma}
The next theorem is somewhat proved in \cite{Q}. Here, for the sake of completeness, we prove the theorem. Note that, our proof gives more information and provides the set of solutions in each case.
\begin{theorem}\label{Square}
	Suppose that $f:\mathbb{Z}_{2^n} \to \mathbb{Z}_{2^n}$ is defined as $f(x)=x^2 \,\,(mod \,\, 2^n)$. Then \\
	\textbf{a)} For three cases $p_2(a)=n $, $p_2(a)=n-1$ with $e_n=0$, and $a=2^{n-2}$ with $e_n=1$  
	$$\vert f^{-1}(a) \vert= 2^\frac{n-1+e_n}{2}.$$  
	\textbf{b)} For two cases $p_2(a)=1 \,\, (mod \,\, 2)$, and $p_2(a)=0\,\,(mod \,\, 2)$ with $0 \leq p_2(a) \leq n-3$ and $o_2(a)\neq 1 \,\, (mod \,\,8)$  
	$$\vert f^{-1}(a) \vert=0.$$  
	\textbf{c)}For the case of $p_2(a)=0 \,\, (mod \,\, 2)$ with $0 \leq p_2(a) \leq n-3$ and $o_2(a)= 1 \,\, (mod \,\, 8)$  
	$$\vert f^{-1}(a) \vert= 2^\frac{p_2(a)+4}{2}.$$  
	\begin{proof}
		\textbf{a)} On one hand, every $a \in \mathbb{Z}_{2^n}$ with $p_2(a) \geq \lceil\frac{n}{2}\rceil$ satisfies $x^2=0 \,\, (mod \,\,2^n)$. So, $\vert f^{-1}(a) \vert$ is at least $2^{n-\lceil\frac{n}{2}\rceil}=2^\frac{n-1+e_n}{2}$. On the other hand, for each $a \in \mathbb{Z}_{2^n}$ with $p_2(a) < \lceil\frac{n}{2}\rceil$, $x^2 \neq 0 \,\, (mod \,\, 2^n)$. Thus, $\vert f^{-1}(a) \vert= 2^\frac{n-1+e_n}{2}$.
		\newline  
		Now, suppose that $n$ is odd and $p_2(a)=n-1$ i.e. $a=2^{n-1}$. Let $x=2^rq$ with odd $q$. We have, 
		$$2^{2r}q^2=2^{n-1}\,\, (mod \,\, 2^n).$$ 
		So $r= \frac{n-1}{2}$, $1 \leq q \leq 2^\frac{n+1}{2}-1$ and $q^2=1 \,\,(mod\,\, 2)$. Thus, only odd $q$'s satisfy the equation $x^2=2^{n-1} \,\, (mod \,\, 2^n)$. Therefore, $\vert f^{-1}(a) \vert= 2^\frac{n-1+e_n}{2}$.
		\newline 
		Now, let $n$ be even and $p_2(a)=n-2$. So $a=s2^{n-2}$, where $s \in \{1,3\}$. If $s=1$, put $x=2^rq$ with odd $q$. Then 
		$$2^{2r}q^2=2^{n-2} \,\, (mod \,\, 2^n).$$ 
		Hence $r= \frac{n-2}{2}$, $1 \leq q \leq 2^\frac{n+2}{2}-1$ and $q^2=1 \,\, (mod \,\, 4)$. So only half of odd $q$'s satisfy the equation $x^2=2^{n-1} \,\, (mod \,\, 2^n)$. Therefore, $\vert f^{-1}(a) \vert= 2^\frac{n-1+e_n}{2}$.
		\newline 
		\textbf{b)} In the proof of the case \textbf{a}, put $s=3$, consider the equation $x^2=2^{n-2} \times 3 \,\, (mod \,\, 2^n)$ and suppose that $x=2^rq$ with odd $q$. Then
		$$2^{2r}q^2=2^{n-2} \times 3 \,\, (mod \,\, 2^n),$$ 
		so $r= \frac{n-2}{2}$ and $q^2=3 \,\, (mod \,\, 4)$. Thus, by Lemma~\ref{aoddmod8}, we have $\vert f^{-1}(a) \vert=0$.
		\newline 
		Now, suppose that $p_2(a)=1 \,\, (mod \,\, 2)$. Since the square of any odd element is odd, so only even elements $x \in \mathbb{Z}_{2^n}$ can satisfy $x^2=a \,\, (mod \,\, 2^n)$. Let $x=2^rq$, $r \neq 0$, and $q$ is odd. Then $p_2(x^2)=2r$ which contradicts $p_2(a)=1 \,\, (mod \,\, 2)$. Therefore $\vert f^{-1}(a) \vert=0$. 
		\newline 
		Now, let $p_2(a)=0 \,\, (mod \,\, 2)$ and $o_2(a) \neq 1 \,\, (mod \,\, 8)$; so $a=2^{2j}t$, where $p_2(a)=2j$ and $t=o_2(a)$. If $x=2^rq$ with odd $q$, then 
		$$2^{2r}q^2=2^{2j}t \,\, (mod \,\, 2^n).$$ 
		Consequently, $r=j$ and $q^2=t \,\, (mod \,\, 2^{n-2j})$. Thus, regarding Lemma~\ref{aoddmod8}, $\vert f^{-1}(a) \vert=0$. 
		\newline 
		\textbf{c)} We use {\rm Theorem 13.3 in~\cite{FR}} to prove this case. Suppose that $p_2(a)=0$ and $a=1 \,\, (mod \,\, 8)$. The algebraic structure $(G,*)$, where $G$ is the subset of odd elements in $\mathbb{Z}_{2^n}$ and $*$ is the operator of multiplication modulo $2^n$ is a group structure. The function $\phi:G \to G$ with $\phi(g)=g*g$ is a group endomorphism on $G$. To compute $|ker(\phi)|$ we must count the number of solutions for the equation $x*x=1_G$. In other words, we must count the number of solutions for the equation $x^2=1 \,\, (mod \,\, 2^n)$. We have
		$$(x-1)(x+1)=0 \,\, (mod \,\, 2^n).$$ 
		Since $x$ is odd, for some $q \in \mathbb{Z}_{2^n}$, $x=2q+1 \,\, (mod \,\, 2^n)$. So
		$$4q(q+1)=0 \,\, (mod \,\, 2^n).$$ 
		Consequently, $q=0$, $q=2^{n-2}$, $q=2^{n-1}$, $q+1=2^{n-2}$, $q+1=2^{n-1}$. Substituting the values of $q$, we have the solutions $x_1=1$, $x_2=2^n-1$, $x_3=2^{n-1}+1$, and $x_4=2^{n-1}-1$. Thus $|ker(\phi)|=4$ and since $|Im(\phi)|=\frac{|G|}{|ker(\phi)|}$, we have, 
		$$|Im(\phi)|=\frac{2^{n-1}}{4}=2^{n-3}.$$ 
		On the other hand, according to Lemma~\ref{aoddmod8} and since the number of elements in $\mathbb{Z}_{2^n}$ in the form of $8q+1$ is equal to $2^{n-3}$ and $|Im(\phi)|=2^{n-3}$, every element in the form of $8q+1$ in $\mathbb{Z}_{2^n}$ is a square. Thus, the equation $x^2=a \,\, (mod \,\, 2^n)$ at least has a solution. Obviously this solution, say $x$, is odd: $x=2y+1$. So we have $(2y+1)^2=8q+1$ or $y^2+y-2q=0$, for some $q$. By Lemma~\ref{aboddceven}, this equation has two solutions $q_1$ and $q_2$. One can check that $q_3=2^n-q_1$ and $q_4=2^n-q_2$ are the two other solutions. Consequently
		$$\vert f^{-1}(a) \vert=|ker(\phi)|=4=2^\frac{p_2(a)+4}{2}.$$  
		Now, suppose that $p_2(a)=0 \,\, (mod \,\, 2),$ and $2 \leq p_2(a) \leq n-3$ and $o_2(a)=1 \,\, (mod \,\, 8)$. In this case, we have $a=2^{2j}t$ with $p_2(a)=2j$ and $t=o_2(a)$. Let $x=2^rq$ with odd $q$. Then
		$$2^{2r}q^2=2^{2j}t \,\, (mod \,\, 2^n).$$ 
		So, $r=j$ and $q^2=t \,\, (mod \,\, 2^{n-2j})$. Regarding Lemma~\ref{aoddmod8} and the proof of case \textbf{b}, this equation has four solutions $q_1,\,q_2,\,q_3,\,q_4$ with $0 \leq q_i \leq 2^{n-2j}-1$. For each of these solutions, we present $2^j$ solutions 
		$$x_{s,t}=2^j \left( s2^{n-2j+1}+q_t \right), \quad 0 \leq s < 2^j, \quad 1 \leq t \leq 4.$$ 
		We have, 
		$$
		\begin{array}{lll}
		x_{s,t}^2=2^{2j} \left( s^22^{2n-4j+2}+q_t^2+2s2^{n-2j+1} \right)\\
		\quad\,\,= s^22^{2n-2j+2}+2^{2j}q_t^2+s2^{n+2}\\
		\quad\,\,=2^{2j}q_t^2 \,\, (mod \,\, 2^n).
		\end{array}
		$$ 
		Regarding $2j \leq n-3$, we have $2n-2j \geq n+3$. Therefore,  
		$$\vert f^{-1}(a) \vert=2^\frac{p_2(a)+4}{2}.$$ 
	\end{proof}
\end{theorem}
\begin{table}\label{table}
	\caption{ The summary of cases of solving equation~(\ref{QuadEqu})}\label{large}  
	\begin{center}
		\begin{tabular}{ccc}
			\hline
			Conditions & Verified in & Number of solutions  \\
			\hline
			$\begin{array}{ll}
			p_2(a)>0,\,p_2(b)>0,\,p_2(c)>0 \\
			t=min \{ p_2(a), \, p_2(b), \, p_2(c) \}
			\end{array}$ & \quad Lemma~\ref{abceven} & \quad $\begin{array}{ll}
			2^t \text{times the number of solutions}\\
			\text{of a corresponding other case} 
			\end{array}$ \\
			$p_2(a)=0,\,p_2(b)=0,\,p_2(c)=0$ & \quad Lemma~\ref{NoSol} & \quad $0$ \\
			$p_2(a)>0,\,p_2(b)>0,\,p_2(c)=0$ & \quad Lemma~\ref{NoSol} & \quad $0$ \\
			$p_2(a)>0,\,p_2(b)=0,\,p_2(c)=0$ & \quad Lemma~\ref{aevenbodd} & \quad $1$ \\
			$p_2(a)>0,\,p_2(b)=0,\,p_2(c)>0$ & \quad Lemma~\ref{aevenbodd} & \quad $1$ \\
			$p_2(a)=0,\,p_2(b)=0,\,p_2(c)>0$ & \quad Lemma~\ref{aboddceven} & \quad $2$ \\
			$\begin{array}{ll}
			p_2(a)=0,\,p_2(b)>0,\,p_2(c)=0 \\
			b=2\acute{b}, \, s=a^{-2}\acute{b}^2-a^{-1}c, \,r=p_2(s)
			\end{array}$ & \quad Corollary~\ref{corollary} & \quad $\begin{array}{ll}
			0 \,\, \text{in some cases}\\
			\text{and}\,\, 2^{{r \over 2}+2} \,\, \text{o.w.} 
			\end{array}$ \\
			$\begin{array}{ll}
			p_2(a)=0,\,p_2(b)>0,\,p_2(c)>0 \\
			b=2\acute{b}, \, s=a^{-2}\acute{b}^2-a^{-1}c, \,r=p_2(s)
			\end{array}$ & \quad Corollary~\ref{corollary} & \quad $\begin{array}{ll}
			0 \,\, \text{in some cases}\\
			\text{and}\,\, 2^{{r \over 2}+2} \,\, \text{o.w.} 
			\end{array}$ \\
			\hline
		\end{tabular}
	\end{center}
\end{table}
Note that Theorem~\ref{Square} gives the set of solutions which are needed in the next corollary.
\begin{corollary}\label{corollary}
	Let $p_2(a)=0$, $p_2(b)>0$, and $b=2\acute{b}$. Put $s=a^{-2}\acute{b}^2-a^{-1}c$, $r=p_2(s)$, and $q=o_2(s)$. If $p_2(r)=0$ or $q\neq1 \,\,(mod\,\, 8)$, then (\ref{QuadEqu}) has no solutions. Otherwise, (\ref{QuadEqu}) has $2^{{r \over 2}+2}$ solutions.
	\begin{proof}
		We have 
		$$ax^2+2\acute{b}x+c=0 \quad (mod \,\, 2^n),$$
		$$x^2+2a^{-1}\acute{b}x+a^{-1}c=0 \quad (mod \,\, 2^n).$$
		Then 
		$$(x+a^{-1}\acute{b})^2=s \quad (mod \,\, 2^n).$$
		Now, by Theorem~\ref{Square}, if $p_2(r)=0$ or $q\neq1 \,\,(mod\,\, 8)$, then (\ref{QuadEqu}) has no solutions and otherwise it has $2^{{r \over 2}+2}$ solutions.
	\end{proof}
\end{corollary}
\begin{remark}
	In Corollary~\ref{corollary}, one should note that if $p_2(c)=2(p_2 (b)-1)$ or $p_2(c)>2(p_2 (b)-1)$ with $p_2 (p_2 (c))=0$, then equation~(\ref{QuadEqu}) has no solutions.
\end{remark}
\begin{remark}
	In Corollary~\ref{corollary}, if $p_2 (c)<2(p_2 (b)-1)$ or $p_2 (c)>2(p_2 (b)-1)$ with $p_2 (p_2 (c))>0$, then we should compute $s \,\, (mod \,\,8)$. The interesting point is that since $a^{-2}=1 \,\, (mod \,\,8)$, it suffices to compute $\acute{s}=\acute{b}^2-a^{-1}c $.
\end{remark}
\section{Conclusion}\label{SecCon}
Quadratic functions have applications in cryptography. In this paper, we study the quadratic equation mod $2^n$. We determine whether this equation has a solution or not and in the case that it has a solution, we give the number of solutions along with the  set of its solutions in $O(n)$ time.
\newline
One of our results is the fact that when the quadratic equation modulo a power of two has a solution, then the number of its solutions is a power of two.

\end{document}